\documentstyle{article}

 \newtheorem{Thm}{Theorem}[section]
 
 \newtheorem{Lem}[Thm]{Lemma}
 \newtheorem{Prop}[Thm]{Proposition}

 \newtheorem{Rem}[Thm]{Remark}
 
 \newcommand{\al}{\alpha}
 \newcommand{\be}{\beta}

 \newcommand{\la}{\lambda}
 \newcommand{\om}{\omega}

 \newcommand{\Ga}{\Gamma}

 \font\msbm=msbm10
 \def\B#1{\hbox{\msbm #1}}
 \def\map#1{\stackrel {#1}\longrightarrow}

 \def\dmap#1{\downarrow #1}

 \newcommand{\n}{\left}
 \newcommand{\r}{\right} 
 \newcommand{\Mo}{(M,\om)}
 \newcommand{\Sym}{Symp(M,\om)}
 \newcommand{\Symo}{Symp_0\Mo}
 
 \newcommand{\pid}{\pi _1(Diff(M))}
 \newcommand{\pis}{\pi _1(\Symo)}
 
 \newcommand{\dk}{\partial _\xi}
 
 \newcommand{\qed}{$\Box$}
 
\begin{document}

  \title{Remarks on the flux groups}
  \author{Jaroslaw K\c{e}dra  
          \thanks {Research supported by the KBN under grants numbers 
          2 PO3A 020 15 and  2 PO3A 02314 .} 
\bigskip \\
           Instytut Matematyczny U. Wr.\\
           pl.Grunwaldzki 2/4\\
           50-384 Wroclaw\\
           Poland\\
           and \\
           Instytut Matematyki U.S.\\
           ul.Wielkopolska 15\\
           70-467 Szczecin\\
           Poland}

  \date{jkedr@math.uni.wroc.pl \\ \bigskip \today}

\maketitle

\begin{abstract}
    Under some topological assumptions, I give a new estimation of the rank of flux
    groups and 
    provide a method of constructions of symplectic aspherical manifolds.
 \end{abstract}

\section{Introduction} \label{S:intro}

\bigskip

Let $\Mo $ be a compact symplectic manifold and $\Sym $ denotes
the group of symplectomorphisms of $\Mo $. Define the
{\bf flux homomorphism }

$$F:\pis \to H^1(M;\B R)$$

$$F(\xi )[a]=\int _{\xi _ta}\om,$$
where $\xi _ta$ denotes the trace of a loop $a$ under the isotopy
$\{\xi _t\}$ representing $\xi $. 
By definition the {\bf flux group } $\Gamma _{\om }$
is the image of the flux homomorphism.

In this paper, we investigate properties the flux groups.
Lalonde, McDuff and Polterovich proved that the rank 
over $\B Z$ of a flux group of $\Mo $ is not greater
than $b_1(M)$ [LMP2, Theorem 2.D].
We extend this result as it is stated
in Theorem A.
The proofs of our results 
(except of the second assertion of Theorem A)
are mainly based on the
work of Gottlieb \cite{go}. Similarly,
Lupton and Oprea used the Gottlieb theory to investigation
of properties of circle actions on, so called, cohomologicaly
symplectic manifolds \cite{lo}.

Let's
briefly look at the contents of the present work.

\bigskip

{\bf Theorem A.}
{\it Let $\Mo $ be a compact symplectic manifold
satisfying one of the following conditions:
   \begin{enumerate}
      \item some Chern number of $\Mo $ is nonzero

      \item $M$ is aspherical and $Z(\pi _1(M))\subset [\pi _1(M),\pi _1 (M)]$

      \item $kerA\subset [\pi _1(M),\pi _1(M)]+torsion$, where 
            $A:\pi _1(M)\to Aut(\pi _*(M))$ is the action of the fundamental
            group on homotopy groups.

   \end{enumerate}
Then 
$\Ga _{\om }\subset ker\n [\cup [\om ]^{n-1}:H^1(M;\B R)\to H^{2n-1}(M;\B R)\r]$
and the rank of $\Ga _{\om }$ over $\B Z$ is not greater than 
dim $ker(\cup [\om ]^{n-1})$.}

\bigskip

We prove this theorem using various methods.
The case (1) follows from the properties of the Wang
sequence associated to a relevant fibration (cf. [LMP2]),
whereas the cases (2) and (3) uses the Gottlieb theory.
The proof of the statement about the rank of flux groups
is a variant of the argument introduced by Lalonde,
McDuff and Polterovich in [LMP2, Theorem 2D].

The above theorem provides an immediate corollary for Lefschetz
manifolds. Recall that a compact symplectic manifold $\Mo $
is called {\bf Lefschetz } if the multiplication
$\cup[\om]^{n-1}:H^1(M;\B R)\to H^{2n-1}(M;\B R)$ is
an isomorphism. Then we have

\bigskip

{\bf Corollary A. } 
{\it Let $\Mo $ satisfies the hypothesis
of Theorem A. If in addition it is Lefschetz then its 
flux group is trivial. }

\bigskip

The next theorem is concerned with aspherical manifolds and
develops the second part of Theorem A.

\bigskip

{\bf Theorem B. } 
{\it Let $\Mo $ be a compact aspherical
symplectic manifold. If either Euler characteristic
of $M$ is nonzero or $\pi _1(M)$ has trivial center
then its flux group $\Ga _{\om }$ is trivial.}

\bigskip

Notice that the above theorem does not depend on the choice
of symplectictic form. The reason is probably that in this
case the groups of symplectomorphisms are 1-connected or
even contractible.

\bigskip

The second section is devoted to construction of aspherical symplectic
manifolds. We use construction of symplectic sums due to Gompf 
\cite{gom} (see also \cite{mw}) and give conditions under which
the resulting manifold is aspherical. At the end, we apply this
construction to give an
example of an aspherical symplectic manifold with zero Chern
numbers and trivial flux group.

\bigskip
\bigskip

{\bf Acknowledgements}

I am grateful to Dusa McDuff for suggesting me the proof
of Theorem A(1), which I previously proved using real
homotopy theory \cite{k} and to Tadek Januszkiewicz
for drawing my attention to the ``Trees'' \cite{se}.

Also, I acknowledge the support of the Foundation for Polish
Science (FNP).

\bigskip

\bigskip

\section{The Proofs} \label{S:tp}
 
\bigskip

Given an element $\xi \in \pis $ there is associated
(up to isomorphism) a symplectic fibration over two sphere 
$\Mo \to P_{\xi }\to S^2$ as follows [LMP2].

$$P_{\xi }=M\times D^2\cup _{\xi _t}M\times D^2,$$ 
where
$D^2$ denotes unit disc and 
$(x,e^{2\pi it})\cong(\xi _t(x),e^{-2\pi it})$.
Fibrations over $S^2$ give raise to the following exact
sequences called the Wang sequences:

$$
...\to H_{k+2}(P_{\xi})\to H_k(M)\stackrel{\dk }{\longrightarrow }
H_{k+1}(M)\stackrel{i_*}{\longrightarrow }H_{k+1}(P_{\xi})\to...
$$

$$
...\to H^k(P_{\xi})\stackrel{i^*}{\longrightarrow }H_k(M)\stackrel{\dk^*}{\longrightarrow }
H_{k-1}(M)\to H_{k+1}(P_{\xi})\to...,
$$

\noindent
where $\dk[a]=[\xi _ta]$ and $\n  <\dk^*\al ,[a]\r>=\n  <\al ,\dk [a]\r>$.
Here $< , >$ denotes the usual pairing between homology and
cohomology. Notice that with the above notation we have that the flux
satisfies the following

$$F(\xi )=\dk ^*[\om ].$$

The crucial property of Wang homomorphism $\dk^*$ is that it is a
derivation \cite{sp}, i.e.

$$\dk^*(\al \cup \be)=
\dk^*(\al )\cup \be+(-1)^{deg(\al )}\al \cup \dk^*(\be ).$$

\bigskip

\begin{Lem}[\cite{dusa84}] \label{L:1}
Let $c_i(M,\om )$ denotes $i^{th}$
 Chern class of 
$M$ with respect to any almost complex structure compatible with $\om $. Then 
$$\dk ^*(c_i(M,\om ))=0,$$
for $\xi \in \pis $.
\end{Lem}

\bigskip

{\bf Proof:} Let $\Mo \to P_{\xi } \to S^2$ be the symplectic fibration 
corresponding
to $\xi $ and consider the subbundle $Vert\subset TP_{\xi }$ consisting of vectors
tangent to the fibers. Since the fibration $P_{\xi }$ is symplectic then 
$Vert$ is symplectic vector bundle which Chern classes restrict to
the Chern classes of $\Mo $: 

$$i^*c_i(Vert)=c_i(i^*Vert)=c_i(TM)=c_i(M).$$ 
Then it follows from the exactness of Wang sequence 
that $\dk ^*(c_i(M,\om ))=0$. 
\qed

\bigskip

Let $C(M)$ denotes the space of continuous maps from
$M$ to $M$ with the open compact topology. 
We define an evaluation $ev_c:C(M)\to M$
by $ev_c(f)=f(x_0)$, $x_0\in M$. Similarly $ev_s$ is
the evaluation defined on $\Sym $. We will use the same
notation for the  maps induced on the fundamental groups.
By $\widetilde {ev_s}$ we denote a homomorphism from
$\pis $ to $H_1(M)=H_1(M;\B Z)/torsion$ which is composition
of $ev_s$ with the obvious maps.

\begin{Lem}[\cite{lmp98}] \label{L:2}
The following diagram is commutative up to
positive constant

\bigskip

$\hskip 2.5cm \pis \buildrel \widetilde {ev_s}\over \longrightarrow  H_1(M) 
\buildrel  PD\over \longrightarrow H^{2n-1}(M)$

$\hskip 3.5cm || \hskip 5cm \downarrow $

$\hskip 2cm \pis \buildrel F\over \longrightarrow H^1(M;\B R) 
\buildrel \cup [\om ]^{n-1}\over \longrightarrow H^{2n-1}(M;\B R).$
  \qed
\bigskip

\end{Lem}

\bigskip

\begin{Rem}
It follows from the above diagram that the flux
conjecture holds for Lefschetz manifolds. Indeed, the image of the
up arrows is discrete in $H^{2n-1}(M;\B R)$ and since
$\cup [\om ]^{n-1}$ is an isomorphism then $\Ga _{\om } $ is also discrete.

\end{Rem}

\bigskip

{\bf Proof of Theorem A.}
First we prove that $\Ga _{\om }\subset ker \cup [\om ]^{n-1}$.

{\bf (1)} Suppose that some Chern number
of $\Mo $ is nonzero. It means that the product of some
Chern classes, say $c_{k_1}\cup...\cup c_{k_m}\neq 0$ in the
top cohomology. It means that $c_{k_1}\cup ...\cup c_{k_m}=\la [\om ]^n$,
for a nonzero $\la \in \B R$. According to the Lemma \ref{L:1} we obtain
that

$$0=\dk ^*(c_{k_1}\cup...\cup c_{k_m})=\dk^*(\la [\om ]^n)=$$
$$=\la n\dk^*[\om ]\cup [\om ]^{n-1}=\la nF(\xi )\cup [\om ]^{n-1},$$
which finishes the first part.

\bigskip

{\bf (2)} Now let $M$ be aspherical and $Z(\pi _1(M))\subset
\n  [\pi _1(M),\pi _1(M)\r] $. Denote by $C_0(M)$
the identity component of $C(M)$. Let's
consider the following diagram

\bigskip
$$
\begin{array}{ccc}
        \Sym  & \map {i}   & C_0(M)\\
 \dmap {ev_s} & \quad      &\dmap {ev_c} \\
       M\qquad & \map {Id}  &M \qquad
\end{array}
$$

\bigskip
\noindent
where the vertical arrows are the evaluations. Gottlieb showed
in [Go, Theorem III.2] that 

$$ev_c:\pi _1(C_0(M))\map {\cong }Z(\pi _1(M)),$$
so passing to the maps on the fundamental groups induced by the above diagram
, we get that 

$$ev_s(\xi )=0$$
 for every $\xi \in \pis .$
Thus we have, due to Lemma \ref{L:2}, that
$F(\xi )\cup [\om ]^{n-1}=0$, which finishes the second part.

{\bf (3)} The proof of this point is the same as previous one since,
once more due to Gottlieb [Go, Theorem I.4], we have that
$ev_c(\pi _1(C_0(M)))\subset ker A$.

Now we have proved that $\Ga _{\om } \subset ker \cup [\om ]^{n-1}$.
For the proof of the second assertion of Theorem A we need the
following lemma (see [LMP2, Theorem 2.A]).

\bigskip

\begin{Lem} \label{L:4}
Let $\om $ and $\om ^{\prime }$ be two symplectic forms on $M$. 
If $\xi \in \pid $ can
be represented by a loops in $\Sym $ and $Symp_0(M,\om ^{\prime })$
then $\dk ^*[\om ]=0$ iff $\dk ^*[\om ^{\prime }]=0$. 
\qed
\end{Lem}

We have to show
that  the rank over $\B Z$
of flux group $\Ga _{\om }$ is not greater than dimension of 
$ker\cup [\om ^{n-1}]$. The proof is in the same vein as the proof
of Theorem 2.D  in 
\cite{lmp99}.

Let $\om _\epsilon $ denotes a rational symplectic form which
is a small perturbation of $\om $. Then
 
$$\hbox {dim} ker\cup [\om ^{n-1}]\geq \hbox {dim} ker\cup[\om _\epsilon ^{n-1}]$$
and 
$c_i(M,\om _\epsilon )=c_i\Mo $. 
Thus $(M,\om _\epsilon )$ satisfies assumptions of Theorem A.

Suppose that there exist elements $\xi _1,...,\xi _k \in \pis$
such that $F(\xi _1),..., F(\xi _k) \in \Ga _{\om }$
are linearly independent over $\B Z$ and
$k>$ dim $ker\cup [\om ^{n-1}]$.
Let 
$\xi _1^{\epsilon },...,\xi _k^{\epsilon }\in \pi _1(Symp(M,\om _\epsilon )$
are represented by  perturbations of loops representing 
$\xi _1,...,\xi _k \in \pis$.
Then the fluxes 
$F_{\epsilon }(\xi _j^{\epsilon })$ 
are rational classes and
are contained in 
$ker\cup[\om _\epsilon ^{n-1}]$, according to already
proved part of Theorem A. 
It follows that some of their
nontrivial combination over $\B Z$ equals zero: 
$\sum _j m_jF_{\epsilon }(\xi _j^{\epsilon })=0$
, $m_j\in \B Z$, $j=1,...,k$.
Due to Lemma \ref{L:4}, we get that
$\sum _j m_j F(\xi _j)=0$ 
which gives the contradiction and completes the proof. 
\qed

\bigskip

 {\bf Proof of Theorem B.} Here we again use the Gottlieb theory,
namely the fact that for a manifold $M$ (not symplectic in general)
 satisfying the hypothesis
of the theorem the space $C_0(M)$ is contractible.
Now suppose that there exists $\xi \in \pis $ with $F(\xi )\neq 0$.
It means that $\dk^*[\om ]\neq 0$, but Wang sequence
depends only on the homotopy type of $P_{\xi }$ and the
latter is homotopy equivalent to the trivial fibration,
due to contractibility of $C_0(M)$. Thus we cannot have
$\dk ^*[\om ]\neq 0$. 
\qed

\bigskip

\bigskip

\section{Symplectic aspherical surgery}

\bigskip

The asphericity of manifolds is a common assumption it the present
paper, so that we shall give a method of construction
of aspherical symplectic manifolds. In fact, the construction
is due to Gompf \cite{gom,mw} and we only give conditions which
ensure the asphericity of the resulting manifold.
 
Let $(M_1,\om _1) $, $(M_2,\om _2) $ and $(S,\om )$ 
be symplectic manifolds,
where $S$ is closed and dim$M_1=$ dim$M_2=$ dim$S-2$.
Let $f_i:S\to M_i$ be symplectic embeddings with normal bundles 
$N_i$, $i=1,2$. Suppose that $c_1(N_1)=-c_1(N_2)$. Then let
$M_1\sharp _SM_2$ denotes the Gompf symplectic sum
of $(M_1,\om _1)$ and $(M_2,\om _2)$ \cite{gom}.

\bigskip

\begin{Prop}[aspherical surgery] 
If in the
above situation
    \begin{enumerate}
          \item $S$ and $M_i-f_i(S)$ are aspherical, $i=1,2.$

          \item the embeddings $f_i$ induce monomorphisms
                $\pi _1(\partial N_i)\to \pi _1(M_i- N_i)=\Ga _i$

\end{enumerate}
Then the symplectic sum $M_1\sharp _SM_2$ is aspherical.
\end{Prop}

\bigskip

{\bf Proof.} The idea of the proof is to show that
the universal covering of  $M_1\sharp _SM_2$ is homotopy
equivalent to a tree. Notice that this fact
is also proved in [Br,proof of Theorem 7.3] by
showing that the higher homotopy groups vanish.

It follows from the second assumption
and van Kampen's theorem
that the fundamental group $\Ga $ of  $M_1\sharp _SM_2$ is
the following amalgamated product $\Ga _1*_{\Ga _0}\Ga _2$,
where $\Ga _0=\pi _1(N_0) $ and $N_0$ is $N_1$ without a the
zero section.
Thus there exists a tree $T$ on which $\Ga $ acts freely
and the fundamental domain of this action is the edge of $T$ 
called the graph
of groups  [Se, Theorem 7]. It follows
from the construction of $T$ that there exists a map $h$
from the universal covering of  $M_1\sharp _SM_2$ to that
tree such that preimage of an egdge is homotopy equivalent
to $\tilde {N_0}$ and preimage of a vertex
is  $\widetilde {M_i-f_i(S)}$ for $i=1$ or 2. Here the tilde
denotes universal covering. Since these preimages
are contractible, due to asphericity, then $h$ is
a homotopy equivalence which ensures the contractibility
of $\widetilde { M_1\sharp _SM_2}$. The latter means that
 $M_1\sharp _SM_2$ is aspherical. 
\qed

\bigskip
The first part of Theorem A may suggest that the flux groups
are nontrivial if the Chern numbers are zero.
In the following example we construct an aspherical symplectic
manifold with zero Chern numbers and trivial flux group.

\bigskip
{\bf Example.}
 Let $M_1=T^2_{\psi }$ be the mapping torus, where $\psi :T^2\to T^2$
is defined by the following matrix

$$
\n [ \begin{array}{cc} 2 & 1 \\ 1 & 1 \\ \end{array} \r ].
$$
In other words, $M_1=\B R^2\times T^2/\sim $, where
$(a+1,b,x,y)\sim (a,b,x,y)$ and $(a,b+1,x,y)\sim (a,b,2x+y,x+y)$.
We endow $M_1$ with the invariant symplectic structure $\om _1$ coming
from the standard symplectic structure on $\B R^4$. This
ensures that Chern classes of $(M_1,\om _1)$ are zero. Moreover, $M_1$ is also
the total space of symplectic $T^2-$bundle over $T^2$, which
means that we can symplectically embed a torus as the fibre.
Let $M_2 = (T^2_{\varphi}, \om _2)$ be another mapping torus
with $T^2$ embedded as a section or a fibre
(in both cases its normal bundle is trivial). Now we perform
a symplectic normal connected sum $M_1\sharp _{T^2}M_2$ \cite{gom,mw}.

Let's now compute fundamental groups that we need.  

$$\pi _1(M_1-N)=(\B Z*\B Z)\times _{\psi }\B Z^2,$$
where $\psi :\B Z*\B Z\to AUT(\B Z^2)$ is defined by
$\psi (a)=Id$ and $\psi (b)=\n [ \begin{array}{cc} 2 & 1 \\ 1 & 1 \\ \end{array}\r ] $,
where $a,b$ are the generators and $N$ denotes the normal
bundle of the fibre.
It is easy to see that the center of the above group
is trivial. The boundary of the normal bundle to the fibre
is three dimensional torus and its inclusion to $M_1-N$ induces
the following monomorphism of the fundamental groups.

$$(k,l,m)\mapsto \n  ((aba^{-1}b^{-1})^k,l,m\r)$$

According to computations of Chern numbers of symplectic sum
in \cite{mw}, we get that the Chern numbers of $M_1\sharp _{T^2}M_2$
are zero because Chern classes of $(M_1,\om _1)$ and $(M_2,\om _2)$
are zero  and Euler characteristic of $T^2$ is zero as well.
Since the center of an amalgamated $G*_KH$ product is equal to
$i^{-1}(Z(G))\cap j^{-1}(Z(H))$, where $i:K\to G,j:K\to H$
are inclusions, then we obtain that $\pi _1(M_1\sharp _{T^2}M_2)$
has trivial center. Due to Theorem B, the flux group of
$M_1\sharp _{T^2}M_2$ is trivial.

Note that the above example is universal in the sense that
a symplectic manifold which is a result of the Gompf sum with 
the manifold $M_1$  
has always trivial center, thus in the aspherical case
its flux is always trivial.

\end{document}